\documentstyle [11pt,amstex]{article}
\def\text#1{{\em #1}}

\def\be{\begin{equation}}
\def\ee{\end{equation}}
\def\bear{\begin{eqnarray}}
\def\eear{\end{eqnarray}}
\def\best{\begin{eqnarray*}}
\def\eest{\end{eqnarray*}}
\def\pf{{\bf Proof}: }

\newtheorem{theorem}{Theorem}[section]
\newtheorem{prop}[theorem]{Proposition}

\newtheorem{lemma}[theorem]{Lemma}

\newtheorem{defn}[theorem]{Definition}

\def\rem{ \addtocounter{theorem}{1}
{\non \bf Remark \arabic{section}.\arabic{theorem} }}

\def\non{\noindent}
\def\pf{\non {\bf Proof. }}
\def\qed{\nopagebreak \hskip .1in { $\Box$ }\penalty10000 %
\hskip\parfillskip \par  }

\def\ra{\rightarrow}

\def\r#1{\right#1}
\def\l#1{\left#1}

\def\ma#1{\mathop {#1} \limits}

\def\Si{\Sigma}

\def\ti{\times}

\def\del{\overline \partial_{j}}

\def\Z{{ \Bbb Z}}

\def\cx{{ \Bbb C}}

\def\ov#1{\overline{#1}}

\def\w{\omega}

\def\M{{\cal M}}

\def\sgn{\mbox{\em sgn}}

\title{\bf The Gromov Invariants of Ruan-Tian and Taubes\vskip.2in}
\author{ Eleny-Nicoleta Ionel\thanks{partially supported by a M.S.R.I.
Postdoctoral Fellowship} \\  M.I.T.\\
 Cambridge, MA  02139 \and Thomas H. Parker\thanks{partially supported by
N.S.F. grant
DMS-9626245}\\ Michigan State University\\ East Lansing, MI   48824}
\date{January 9, 1997}

\addtocounter{section}{0}
\begin{document}

\maketitle

\vskip.4in

C. Taubes has recently defined Gromov invariants for symplectic four-manifolds
and related them to the
Seiberg-Witten invariants ([T1], [T2]).  Independently, Y. Ruan and G. Tian
defined symplectic invariants based on ideas of Witten ([RT]).  While
similar in spirit, these two sets of
invariants are quite different in their details.

In this note we show that Taubes' Gromov invariants are equal to certain
combinations of Ruan-Tian invariants
(Theorem \ref{4.mainthm}). This link allows us to generalize Taubes'
invariants.  For each closed symplectic
four-manifold, we define a sequence of symplectic invariants $Gr_{\delta}$,
$\delta=0,1,2\dots$.  The first of these, $Gr_0$,
generates Taubes' invariants, which count embedded $J$-holomorphic curves.
The new invariants $Gr_{\delta}$ count
immersed curves with
$\delta$ double points.

In particular, these results give an independent proof that Taubes'
invariants are well-defined.
Combined with Taubes' Theorem [T1], they also show that, for symplectic
4-manifolds with $b^+>1$, some of the
Ruan-Tian symplectic invariants agree with the Seiberg-Witten invariants.

\medskip
\section{Gromov Invariants}

 Fix a closed symplectic four-manifold $(X,\w)$.  Following the ideas of
Gromov and Donaldson, one can define
symplectic invariants by introducing an almost complex structure $J$  and
counting (with orientation) the number
of $J$-holomorphic curves on $X$ satisfying certain constraints.
Unfortunately, technical difficulties make it
necessary  to modify the straightforward count in order to obtain an
invariant.  In this section we review the
general construction and describe how the technicalities have led to two
types of Gromov invariants.

Given $(X,\w)$, one can always choose an  almost complex
structure $J$ tamed by $\w$, i.e. with $\w(Z,JZ)>0$ for all tangent vectors
$Z$.  A map $f:\Si \to
X$ from a topological surface $\Si$ is called  $J$-holomorphic if there is
a complex structure $j$ on $\Si$ such
that
\begin{equation}
\overline{\partial}_Jf=0
\label{1.holomorphicmapeq}
\end{equation}
where $\overline{\partial}_Jf =\frac12(df \circ j - J \circ df)$.  The
image of such a map is a  $J$-holomorphic
curve. Conversely, each  immersed $J$-holomorphic curve is uniquely
specified by the equivalence class of a
$J$-holomorphic pair $(f,j)$ under the action of the group of
diffeomorphisms of $\Si$.  These equivalence
classes $[(f,j)]$ form moduli spaces
$$
{\cal M}_{{A,g}}
$$
labeled by the genus $g$ of $\Si$ and the class $A\in H_2(X)$ of the image
(and implicitly depending on $J$).
The formal tangent space to ${\cal M}_{{A,g}}$ at $[(f,j)]$ can be identified
with the kernel of the operator
\begin{equation}
D_{f,j}:\Gamma(f^*N)\to \Omega^{0,1}(f^*N)
\label{1.kernelD}
\end{equation}
obtained by  linearizing (\ref{1.holomorphicmapeq}) and restricting to the
normal bundle $N$ along the image of
$f$.  The Riemann-Roch Theorem  shows that
$$
\mbox{dim}\ {\cal M}_{{A,g}} =2[g-1-\kappa\cdot A]
\label{1.dimM}
$$
where $\kappa$ is the canonical class of $(X,J)$. We can elaborate on this
construction by marking $d$
points $x_i$ on $\Si$.  The equivalence classes $[(f,j,x_1,\dots ,x_d)]$ of
marked $J$-holomorphic curves then
form a moduli space ${\cal M}_{{A,g,d}}$ of dimension $2[g-1-\kappa\cdot
A+d]$, and the evaluations $x_i\mapsto
f(x_i)$ define a map
$$
ev:{\cal M}_{{A,g,d}} \to X^d=X\ti\dots\ti X
$$

 The marked points enable us to restrict attention to curves satisfying
constraints. For our purposes it is almost always
enough to consider point constraints.   Thus we pick
$$
d=d_{A,g}=g-1-\kappa\cdot A
$$
generic points $p_i$ in $X$ and consider the constrained moduli space
$$
{\cal M}'_{{A,g}}=ev^{-1}(p_1,\dots, p_d)
$$
of all $J$-holomorphic curves that pass through the (ordered) points $p_i$.
For generic $J$ and $\{p_i\}$, this
constrained moduli space is zero-dimensional and its formal tangent space
at $C=[(f,j,x_1,\dots ,x_d)]$ is the
kernel of the restriction
$D_C$ of (\ref{1.kernelD}) to the subspace of $\Gamma(f^*N)$ that vanishes
at each marked point. Each curve $C\in
{\cal M}'_{{A,g}}$ therefore has a sign given by
$(-1)^{\em SF}$ where SF denotes the spectral flow from $D_C$ to any
complex operator $\overline{\partial}_C$
which is a compact perturbation of $D_C$.  Counting the points in ${\cal
M}'_{{A,g}}$ with sign gives a ``Gromov
invariant''
\begin{equation}
Gr_{A,g}(p^d)=\sum_{C\in{\cal M}'_{{A,g}}}{\mbox{sgn}\ C} =
\sum_{C\in{\cal M}'_{{A,g}}}\
(-1)^{\em{SF}\,(D_C)}.
\label{1.SF}
\end{equation}

One then tries to mimic  Donaldson's cobordism arguments to show that
$Gr_{A,g}$ is independent of
$J$ and $\{p_i\}$, and hence defines a symplectic invariant. This involves
considerable
analysis, and along the way  one encounters a major technical difficulty
---  $\M_{A,g}$ may not be a manifold at
the multiply-covered maps.   There currently exist two distinct ways of
overcoming this difficulty.

\begin{enumerate}
\item   Taubes restricts  $g$ to be the genus expected for embedded curves and
counts embedded, not necessarily connected, $J$-holomorphic curves, dealing
with the complications associated
with multiply-covered curves.  In the end he obtains `Gromov-Taubes'
invariants that we will denote by $GT_0(A)$.

\item   Ruan-Tian [RT] observed that the
difficulties with  multiply-covered maps can be overcome by replacing
(\ref{1.holomorphicmapeq}) by the
inhomogeneous equation
$$
\overline{\partial}_Jf=\nu
\label{1.pertholomorphicmapeq}
$$
where  $\nu$ is an appropriate perturbation term.    We will denote the
resulting symplectic invariants by
$RT(A,d)$.
\end{enumerate}

The next two sections give some details about  these two sets of invariants  and
describe generating functions involving them.

\medskip

\medskip
\section{The Taubes Series}

The details of Taubes' construction are interesting and surprisingly
subtle.  Given  $A\in
H_2(X,\Z)$, Taubes fixes  the genus
to be
$$
g_A=1+\frac12(A\cdot A +\kappa\cdot A).
$$
The moduli space of such curves has   $\mbox{dim}\, \M_A =A\cdot
A-\kappa\cdot A$, so we constrain by  $d_A=\frac12
(A\cdot A-\kappa\cdot A)$ points. The adjunction formula implies that each
constrained  curve is embedded unless
$A$ lies in the set
$$
T=\{\; A\in H_2(X,\Z)\;|\; A^2=0 \ \ \mbox{and}\ \ \kappa\cdot A=0\;\},
$$
in which case the curve is a multiple cover of an embedded torus and
$d_A=0$. Similarly, each constrained  curve
in
$$
{\cal E}=\{\; A\in H_2(X,\Z)\;|\; A^2=-1\;\},
$$
is an embedded ``exceptional'' sphere.

More generally, for each class $A$ and $d\geq 0$ we get a  count
of {\it connected} curves through $d$ generic points
\best
Gr(A,d)
\label{2.AnotinS}
\eest
defined by (\ref{1.SF}) with $g=d+1+\kappa\cdot A$. Note that by the
adjunction formula
\begin{equation}
d_A-d=g_A-g=\delta \geq 0,
\label{2.defDelta}
\end{equation}
so $0\le d\le d_A$ with $d_A=0$ for $A\in {\cal E}\cup T$.  Geometrically,
$\delta$ is the
number of double points on a generic immersed $A$-curve.

 Taubes observed that for $A\in T$, $Gr(A,0)$  depends on $J$, as follows.
For an embedded torus $C$, let
$L_i,\ i=1,2,3$ be the three non-trivial real line bundles over $C$.
Twisting the linearization $D_C$ by $L_i$
gives operators
$$
D_i:\Gamma(f^*N\otimes L_i)\to \Omega^{0,1}(f^*N\otimes L_i).
$$
The space of almost complex structures is divided into chambers by the
codimension one
``walls''  consisting of those $J$ for which there is a $J$-holomorphic
curve with either $D_C$ or one of the
$D_i$ not invertible.  The value of $Gr(A,0)$ changes as $J$ crosses a wall.

Within a chamber, there are four types of $J$-holomorphic tori,  labeled by
the number $k=0,1,2,3$  of the $D_i$
whose sign (determined by the spectral flow) is negative.   Thus for
generic $J$, the moduli space of
$J$-holomorphic
$A$-curves  is the disjoint union of four zero-dimensional moduli spaces
${\cal M}_{A,k}$.  Counting with sign
gives four ``Taubes numbers''
\begin{equation}
\tau(A,k)=  \sum_{C\in{\cal M}_{A,k}}\ \mbox{sgn}\ C.
\label{2.AinS}
\end{equation}
Taubes  derived  wall-crossing formulas and showed that a certain
combination of the $\tau(A,k)$ is independent of
$J$.

The right combination is best described by assembling the  counts
(\ref{2.AnotinS}) and (\ref{2.AinS})  into a
single quantity associated with $X$.  For that purpose, we introduce formal
symbols $t_A$ for $A\in H_2(X;\Z)$
with relations
$t_{A+B}=t_At_B$ and specify three ``generating functions'' $e(t),f(t)$ and
$g(t)$.   From $f$ we construct
functions $f_k$ corresponding to the four types of curves by setting
\bear\label{T gen fc}
f_0=f,\qquad
f_1(t)={f(t)\over f(t^2)},\qquad
f_2(t)={f(t)f(t^4)\over f^2(t^2)},\qquad
f_3(t)={f(t)f(t^4)\over f^3(t^2)}.
\eear
in accordance to the wall crossing formulas  in [T2].  We will also use
another variable $s$ to keep track of
the number of double points.

\begin{defn}  The {\em Taubes Series} of $(X,\w)$ with generating functions
$e, f$
and $g$ is the formal power series in the variables $t_A$ and $s$ defined by
\bear
GT_X(t,s)=\prod_{E\in {\cal E}}e(t_E)^{Gr(A,0)} \cdot
\prod_{A\notin T\cup {\cal E}}\prod_{d=0}^{d_A}
g\l(t_A{s^{d}\over d! }\r)^{Gr(A,d)}
\cdot \prod_{A\in T}\ma\prod_{ k=0}^3 f_k(t_A)^{\tau(A,k)}
\label{egfgeneratingfnc}
\eear
with  the $f_k$ given by (\ref{T gen fc}).
\end{defn}
We then get a sequence of maps $GT_{\delta}: H_2(X;\Z)\to \Z$ by expanding
(\ref{egfgeneratingfnc}) as a power series in
$s$:
\bear
GT(t,s)\ =\ \sum_{A} \sum_{\delta=d_A-d}GT_{\delta}(A)\ t_A\, \frac{s^d}{d!}
\label{egfgeneratingfnc2}
\eear
where we have labeled the coefficients by $\delta=d_A-d$ rather than $d$.

\begin{prop} With the choice
\bear
e(t)=1+t,\qquad  f(t)=\frac{1}{1-t},\quad \mbox{and}\qquad  g(t)=e^t,
\label{2.taubesgeneratingfncs}
\eear
the degree zero component $GT_{0}$ in (\ref{egfgeneratingfnc2}) is the
Gromov invariant defined by Taubes in
[T2].
\end{prop}

\pf The coefficient $GT_{\delta}(A)$ of $t_A s^d/d!$ in
(\ref{egfgeneratingfnc}) is a sum of coefficients, one for
each  product of monomials $(t_{A_i}s^{d_i})^{n_i}$ with $d=\sum n_id_i$
and $A=\sum n_iA_i$, where the $A_i$ are
distinct homology classes, $n_i\geq 0$, and $n_i=1$ for all  $A_i\in {\cal
E}$ (because the generating function
is  $e(t)=1+t$).  Given such a decomposition, we can expand $\delta=d_A-d
=d_A-\sum n_id_i$ by writing
$\delta_i=d_{A_i}-d_i\ge 0$ as in (\ref{2.defDelta}) and using the
definition of $d_A$.  This gives
\best
\delta  & = & \frac12\left[(\ma\sum n_i A_i)^2-\ma\sum n_i A_i^2\right]
+\sum n_i \delta_i \\
 & = & \ma\sum \frac12 n_i(n_i-1)\ A_i^2 +\sum_{i<j} n_i n_j A_i A_j
+\sum n_i\delta_i
\eest
Each of the terms in this sum are nonnegative since (a) $A_i^2\ge 0$ for
$A_i\notin {\cal E}$ and
$n_i=1$ for  $A_i\in {\cal E}$, and (b) $A_i\cdot A_j\ge 0$ for $i\ne j$
because the $A_i$ are distinct.
Consequently, the only monomials that contribute to the $\delta=0$ term
are those corresponding to
decompositions of $A$ and $d$ with
\smallskip

\ \hskip1in (a) $n_i=1$ unless $A_i^2=0$,

\ \hskip1in  (b) $A_i\cdot A_j=0$ for all $i\ne j$,

\ \hskip1in  (c) $d_i=d_{A_i}$.

\smallskip

\non Let ${\cal S}={\cal S}(A)$ be the set of such decompositions.  For
each $y=\{(n_i,A_i)\}$ in ${\cal S}$, let $y'$ be the set of those
 $(n_i,A_i)\in y$ with $A_i\notin T$, let $y''$ be the set of those
 $(n_i,A_i)\in y$ with $A_i$ primitive and $A_i\in T$,
 and let $t_{y'}$ and $t_{y''}$ be the corresponding monomials.
Putting the functions (\ref{2.taubesgeneratingfncs}) into
(\ref{egfgeneratingfnc}), one sees that the
coefficient of $t_A s^{d_A}/d_A!$ has the form
\bear
\label{2.T1}
Gr_0(A)\ =\ \sum_{y\in{\cal S}} R(y') Q(y'').
\eear
Here $R(y')$ is the coefficient of $t_{y'}s^d/d!$ in
$$
\prod_{A_i\notin T}
\left[\mbox{exp}\left(t_{A_i}\frac{s^{d_i}}{d_i!}\right)\right]^{Gr(A_i,d
_i)}
$$
(after noting that $t_{y'}$ is at most linear in $t_{A_i}$ for each $A_i\in
{\cal E}$ and $\mbox{exp}\,t=1+t+O(t^2)$), so
\bear
\label{2.T2}
R(y')\ =\ d!\,\prod_{(n_i,A_i)\in y'}  \frac{Gr(A_i,d_{i})^{n_i}}
{n_i!\,(d_i!)^{n_i}}.
\eear
Similarly,  $Q(y'')$ is the coefficient of $t_{y''}$ in
$$
\prod_{A_i\in T}\ \prod_{ k=0}^3 f_k(t_{A_i})^{\tau(A_i,k)}=
\prod_{A_i\in T \atop primitive}\ \prod_{q=1}^{\infty}\;
\ma\prod_{ k=0}^3 f_k(t_{A_i}^q)^{\tau(qA_i,k)}.
$$
Then
\bear
\label{2.T3}
Q(y'')\ =\ \prod_{(n_i,A_i)\in y''}\ Q(n_i,A_i)
\eear
where   $Q(n,A)$ is the coefficient of $t_{A}^n$ in
$$
\ma\prod_{q=1}^{\infty}\ma\prod_{ k=0}^3 f_k(t_{qA})^{\tau(qA,k)}.
$$
 For each  embedded, holomorphic torus $C$, let $f_C$ denote the function
$f_k$ (resp. $1/f_k$) when $C$ is of type $k$ and has positive (resp.
negative) sign.
Expanding $f_{C}(t)=\sum_{m} r(C,m) t^m$, we have
\bear
\label{2.T4}
Q(n,A)\ =\ \sum_{{\cal D}}\prod r(C_j,m_j),
\eear
where  ${\cal D}$ is the set of all pairs $(m_j,C_j)$ of
$J$-holomorphic curves $C_j$ and
multiplicities $m_j$ with $[C_j]=q_jA$  and $\sum m_j q_j=n$.
Together, (\ref{2.T1}) --  (\ref{2.T4}) exactly agree with the invariant
defined by Taubes ([T2] section 5d). \qed

\bigskip

\begin{rem}  Taubes chooses the functions (\ref{2.taubesgeneratingfncs}) to
make his invariants agree with the
Seiberg-Witten invariants.
\end{rem}

\medskip

The numbers $Gr_{\delta}(A)$ defined by (\ref{egfgeneratingfnc}) and
(\ref{egfgeneratingfnc2}) count the
$J$-holomorphic $A$-curves (of any genus and any number of components) with
$\delta$ double points, and thus generalize Taubes'
count of embedded curves. Below, we will verify that  the
$Gr_{\delta}(A)$ are symplectic
invariants by  relating the Taubes Series to Ruan-Tian invariants.

\medskip

\medskip
\setcounter{equation}{0}
\section{The RT Series}

Ruan and Tian [RT] define symplectic invariants
$RT_{A,g,d}(\alpha_1,\dots,\alpha_d)$ by taking the moduli space
$\M_{A,g,d}$ of {\em connected, perturbed} holomorphic  $A$-curves with
genus $g$ and $d$ marked points,  restricting
to the subset
$\M'_{A,g,d}$ where the marked points lie on fixed constraint surfaces
representing the $\alpha_i\in H_*(X)$, and counting with orientation
(assuming  $\M'_{A,g,d}$ is
zero-dimensional). In particular,
when the $\alpha_i$ are all points and $g=d+1+\kappa\cdot A$ we get invariants
\begin{equation}
RT(A,d)\ =\ RT_{A,d+1+\kappa\cdot A,d}(p^d).
\label{3.1}
\end{equation}
This section describes how to assemble these invariants into a series
analogous to
(\ref{egfgeneratingfnc}).

First we must deal with a technical problem.  In [RT], the invariants
$RT_{A,g,d}$ are defined  only for the
``stable range'' $2g+d\geq 3$. This leaves $RT(A,d)$ undefined for two
types of curves: tori with no marked
points, which occur when $d=\kappa\cdot A=0$,  and  spheres
with fewer than three marked points, which occur when  $d=0,1,2$ and
$d+1=-\kappa\cdot A$. But we can extend
definition (\ref{3.1}) to these cases by imposing  additional
``constraints'' which are automatically satisfied.  For this,  choose a
class $\beta\in H_{2}(X)$  with
$A\cdot \beta \neq 0$ and set
\begin{equation}
RT(A,0) \ =\   \frac{1}{A\cdot \beta}\, RT_{A,1,1}(\beta) \qquad \mbox{if}\
\kappa\cdot A =0
\label{3.2}
\end{equation}
and
$$
RT(A,d) \ =\   \frac{1}{(A\cdot
\beta)^{3-d}}\,RT_{A,0,3}(p^d\beta^{3-d})\qquad  \mbox{if }\ d=\kappa\cdot
A-1=0,1,2.
$$
Thus defined, these invariants  count
perturbed holomorphic curves.  For example, when $\kappa\cdot A =0$ each
genus one curve $C$
(without marked points) representing
$A$ is a map $f:T^2\to X$, well-defined up  the automorphisms of $T^2$ with
the induced
complex structure.  Fix a point $p\in T^2$ and represent $\beta$ by a cycle
in general position.  Then $C\cap
\beta$ consists of  $A\cdot \beta$ distinct points.  Hence $C$ is the image
of  exactly $A\cdot \beta$ maps $f:T^2\to X$ with $f(p)\in\beta$ and these
are counted by $RT_{A,1,1}(\beta)$.

\medskip

Now fix  a generating function $F_A$ for each class $A$  and assign a
factor $F_A(t_A)$ to each curve that
contributes $+1$ to the count $RT(A,d)$, and a factor $1/F_A(t_A)$ to each
curve that
contributes $-1$.  Taking the product gives  a   series in the variables $t_A$
$$
\ma\prod_{A\in H_2(X)}F_A(t_A)^{RT(A,d)}
$$
which is an invariant of the deformation class of the symplectic structure
of $(X,\w)$. As with the Taubes Series,
different choices of the  $F_A$ give different  series, but all  encode the
same data. We will choose three
generating functions  and form a series resembling (\ref{egfgeneratingfnc}).

\begin{defn} The {\em Ruan-Tian Series} of $(X,\w)$ defined by $e(t)$,
$F(t)$ and $g(t)$ is
\bear\label{3.defgr}
RT_X(t,s)\ =\ \prod_{E\in {\cal E}}e(t_E)^{RT(A,0)} \cdot
\prod_{A\notin T\cup {\cal E}} g\l(t_A{s^{d}\over d! }\r)^{RT(A,d)}\ \cdot\
\ma\prod_{A\in T}F(t_A)^{RT(A,0)}
\eear
Expanding in power series as in (\ref{egfgeneratingfnc2}) gives invariants
$RT_{\delta}:H_2(X;\Z)\to\Z$.
\end{defn}

To make this more concrete, we could take $e(t)$, $F(t)$ and $g(t)$ to be
the specific functions given in
(\ref{2.taubesgeneratingfncs}).  That choice, however, overcounts tori with
self-intersection zero.   It
turns out that the formulas are simpler if $F$ satisfies
\bear\label{prodF=t}
\prod_{k=1}^\infty F(t^k)=e^t.
\eear
Thus it is appropriate to make the more awkward-looking choice
\bear\label{MoebiusF}
e(t)=1+t,\qquad F(t)=\exp\l(\ma\sum_{m= 1}^\infty \mu(m)t^m \r) ,\qquad
\mbox{and}\qquad  g(t)=e^t,
\eear
where $\mu$ is the M\"{o}bius function.  (The M\"{o}bius function is
defined by $\mu(1)=1$, $\mu(m)=(-1)^k$ if
$m$ is a product of $k$ distinct primes, and $\mu=0$ otherwise.)  One can
then verify (\ref{prodF=t}) by writing
$\ell=mk$ and using the basic fact that
$$
 \ma\sum_{m|\ell}\mu(m)\ =\ \left\{
\begin{array}{ll}
1\qquad & \mbox{if }\ell=1,\\
0 & \mbox{otherwise.}
\end{array}
\right.
$$
 We will see next how the generating functions (\ref{MoebiusF}) lead back
to the Taubes Series and the
Seiberg-Witten invariants.

\medskip

\medskip
\setcounter{equation}{0}
\section{Equivalence of the Invariants}

In this section we will prove that the Taubes and Ruan-Tian Series are equal
for any closed symplectic four-manifold. The proof is straightforward for
classes $A\notin T$, but for the
toroidal classes $A\in T$ it requires some combinatorics.

\smallskip

For classes $A\notin T$, the moduli space of
$J$-holomorphic curves of genus $g_A$  passing through $d$ points
contains no multiply covered curves for generic $J$ (cf. [R], [T2]).
Consequently, the
moduli space of such curves is smooth  and the linearized operator has no
cokernel.  The Implicit Function Theorem then implies that  each of these
curves (but none of their multiple
covers) can be uniquely perturbed to a solution of the equation $\del
f=\nu$ for small $\nu$.  Thus
\bear
Gr(A,d)=RT(A,d)\qquad \mbox{ for } \quad A\notin  T,
\label{4.1}
\eear
so the first two factors in the products (\ref{egfgeneratingfnc}) and
(\ref{3.defgr}) are equal.

The computations for $A\in T$ are more complicated because multiple covers
{\em do} contribute. In this case,
the moduli space
$\M_A$ of $J$-holomorphic, connected, embedded $A$-curves is finite for
generic $J$, and each curve $C\in \M_A$ is a torus.  The last part
of the Gromov series
(\ref{3.defgr})   has the form
$$
Gr^T\ =\ \prod_{A\in T} \ \prod_{C\in \M_A} \phi_C(t_A)
$$
for some function $\phi_C$ that we must determine.

To do that, we fix one  torus $C\in\M_A$ defined by an embedding
$(T^2,x_0,j_0)\ra X$ and regard the
domain $(T^2,x_0,j_0)$ as the quotient of the complex plane by the lattice
$$
\Lambda_0=\Z\oplus \tau\Z.
$$
Curves $C'$ which are $m$-fold covers of $C$ are given  by pairs
$(\psi,j)$ where  $\psi:(T^2,x_0,j)\ra (T^2,x_0,j_0)$ is an $m$-fold  cover map;
these are classified  (up to
diffeomorphisms of the domain) by  index $m$ sublattices $\Lambda \subset
\Lambda_0$.  Let ${\cal L}_m$ be the set
of all such lattices.

For generic $J$ the linearized operator has zero cokernel (it is invertible
with index zero).  Hence
each  $m$-fold cover  can be uniquely
perturbed to a solution of  $\del f=\nu$, which contributes to
$RT_{mA,1,1}$.  The total contribution of the multiple covers of $C$ to
$RT_{mA,1,1}$ is
$$
\ma\sum_{\Lambda\in{\cal L}_m} \mbox{sgn }\Lambda
\label{sumofsgns}
$$
where  $\mbox{sgn }\Lambda$ is the sign of the multiple cover $C'$
described by $\Lambda$.  Thus, after  stabilizing as in
(\ref{3.2}),
\bear
\phi_C(t_A)\ =\ \prod_{m=1}^{\infty}F(t_A^m)^{\frac{1}{m}
\ma\sum_{\Lambda\in{\cal L}_m} \mbox{sgn }\Lambda}
\label{phiC}
\eear

 To proceed, we must
determine $\mbox{sgn }\Lambda$ using the orientation prescribed by
Ruan-Tian.  As in Section 1, this is given by the the spectral
flow of the linearization $D_C$ (the exposition in
[RT] is obscure, but this is clearly the orientation
that the authors intended to specify).  This sign is independent of $\nu$
for small $\nu$, so we can assume that $\nu=0$ in the subsequent
 calculations.

\begin{lemma}
The sign of a curve $C'=\cx/\Lambda$ is
\begin{equation}
\sgn \, \Lambda\ =\ \sgn \,D_0\,\ma\prod \sgn\,D_i
\label{sgnLambda}
\end{equation}
where the product is over all $i=1,2,3$ such that $\Lambda_0$ is a sublattice
of $\Lambda_i$ with $\Lambda_i$  defined by (\ref{lambdai}).
\end{lemma}

\pf  Looking at the explicit formula for $D_{C'}$ [T2],
one sees that $D_{C'}$ is the pullback of $D_{C}$ (it depends only on the
1-jet of $J$ along $C$). Fix a complex
operator
$\ov\partial$ on $C$, choose a path from
$\ov\partial$ to $D_{C}$, and let $D_t$ be the lifted path of operators on
$C'$;  each $D_t$ is  invariant
under deck transformations.  As in [T2], we can assume that $\mbox{ker}\,
D_t=\{0\}$ except at
finitely many values of $t=t_k$, where $\mbox{ker}\, D_t$ is one-dimensional.

The translations of $\cx$ by 1 and $\tau$ respectively
induce  deck transformations $\tau_1$ and $\tau_2$ of $C'\to C$; these
generate the abelian group
$G=\Lambda_0/\Lambda$ of all deck transformations.  At each $t=t_k$,
$\mbox{ker}\, D_t$ is a one-dimensional representation $\rho_i$ of $G$, so is
one of four possibilities:
$$
\left\{\begin{array}{l}\rho_0\tau_1(\xi)=\xi\\
\rho_0\tau_2(\xi)=\xi\end{array}\right.\qquad
\left\{\begin{array}{l}\rho_1\tau_1(\xi)=-\xi\\
\rho_1\tau_2(\xi)=\xi\end{array}\right.\qquad
\left\{\begin{array}{l}\rho_2\tau_1(\xi)=\xi\\
\rho_2\tau_2(\xi)=-\xi\end{array}\right.\qquad
\left\{\begin{array}{l}\rho_3\tau_1(\xi)=-\xi\\
\rho_3\tau_2(\xi)=-\xi\end{array}\right.
$$
where $\xi$ is a generator of the kernel.  Call these kernels of type
0,\,1,\,2 and 3 respectively. Then
$$
SF=\sum_{i=0}^3 SF_i
$$
where $\mbox{SF}_i$ is the  number of $t_k$ of type $i$ (counted with
orientation), and
\begin{equation}
\mbox{sgn } \Lambda\ =\ (-1)^{\mbox{SF}}\ =\ \prod (-1)^{\mbox{SF}_i}.
\label{sgnLambda2}
\end{equation}
Note that each $\xi$ of type 0 descends to a section of $\mbox{ker
}\,D_t$ on $\cx/\Lambda_0$.  In fact, this
is a one-to-one correspondence, so $SF_0$ is the spectral flow of the path
$D_t$ on the base curve $C$ and
$(-1)^{\mbox{SF}}$ is the sign of $D_0$. The remaining representations
determine three index two sublattices
\bear
\Lambda_i=\mbox{ker}\, \rho_i
\label{lambdai}
\eear
of $\Lambda$.  Thinking of $\xi$ as a $\Lambda$-invariant section on $\cx$,
one sees that  (a) a type $i$ kernel
cannot appear unless $\Lambda\subset\Lambda_i$, and (b) if
$\Lambda\subset\Lambda_i$ then $\xi$  descends to an
element of $\mbox{ker}\, D_t$  over the double cover $\cx/\Lambda_i$.  Thus
$\mbox{SF}_i$ vanishes if
$\Lambda$ is not a subset of $\Lambda_i$, and when
$\Lambda\subset\Lambda_i$ $\mbox{SF}_i$ coincides with the
spectral flow of Taubes' operator $D_i$.  Then (\ref{sgnLambda2}) is the
same as (\ref{sgnLambda}).
\qed

\bigskip

\rem A set of representatives of the lattices in ${\cal L}_m$ is
\begin{equation}
\Lambda=a\Z+(b\tau+p)\Z  \qquad\mbox{where}\ \ m=ab,\ \ p=0,\dots,a-1.
\label{lattices}
\end{equation}
The group of deck transformations is $G\cong \Z_{a}\ti\Z_{b}$, and the
three the lattices (\ref{lambdai}) are
$$
\Lambda_1=\Z+2\tau\Z, \qquad  \Lambda_2=2\Z+\tau\Z, \qquad
\Lambda_3=2\Z+(1+\tau)\Z.
$$

\medskip

Now fix $m$ and separate the set of lattices  ${\cal L}$ into:
\best
{\cal L}^0\ &=&\{ \Lambda\in{\cal L}\;|\; \Lambda
\mbox{ is contained in none of the lattices }\Lambda_1,\Lambda_2,\Lambda_3\;\}\\
{\cal L}^i\ &=&\{ \Lambda\in{\cal L}\;|\; \Lambda
\mbox{ is contained in $\Lambda_k$ only for $k=i$}\;\}\\
{\cal L}^{123}&=&\{ \Lambda\in{\cal L}\;|\; \Lambda
\mbox{ is contained in }\Lambda_1,\Lambda_2\ \mbox{and}\ \Lambda_3\;\}.
\eest
Note that if $\Lambda$ is contained in  two of the $\Lambda_i$ then it is
contained in the third.  Thus the
above sets constitute a partition
$$
{\cal L}= {\cal L}^0 \cup {\cal L}^1 \cup {\cal L}^2 \cup {\cal L}^3 \cup
{\cal L}^{123}.
$$
Furthermore, there are automorphisms of $\Lambda_0$ that interchange the
lattices $\Lambda_i, i=1,2,3$, so the
sets ${\cal L}^1$, ${\cal L}^2$, and ${\cal L}^3$ have the same
cardinality.   Hence from  (\ref{sgnLambda}) we have
\bear
\sum_{\Lambda\in{\cal L}_m}  \mbox{sgn }\Lambda\ =\ \mbox{sgn}\; D_0\left\{
A+B\;\sum_{i=1}^3 \mbox{sgn}\;
D_i+C\prod_{i=1}^3
\mbox{sgn}\; D_i \right\}
\label{RTm3}
\eear
where $A=|{\cal L}^0|$ is the number of elements of ${\cal L}^0$,
$B=|{\cal L}^1|$, and $C=|{\cal L}^{123}|$.

\medskip

\begin{lemma}
\label{abclemma}
 Set $\sigma(m)=\ma\sum_{a|m}a$ if $m$ is a
positive integer, and $\sigma=0$ otherwise.  Then
$$
A+3B+C=\sigma(m), \qquad B+C=\sigma(m/2), \qquad C=\sigma(m/4).
$$
\end{lemma}
\pf Using the representatives (\ref{lattices}) of ${\cal L}$, we have
\best
A+3B+C=|{\cal L}|=\sum_{m=ab}\;\sum_{p=0}^{a-1}
1=\ma\sum_{a|m}a=\sigma(m).
\eest
Next, $B+C$ is the number of lattices $\Lambda\in{\cal L}$ which contain
$\Lambda_1$.  These are the lattices
(\ref{lattices}) with $b=2\beta$ even, so
\best
B+C=\sum_{m=a\cdot2\beta}\;\sum_{p=0}^{a-1}1=
\sum_{a|{m\over 2}} a=\sigma\l({m\over 2}\r).
\eest
Finally, ${\cal L}^{123}$ is the set of all lattices $\Lambda$ such that
$a, b$ and $p$ are all even. Writing
$a=2\alpha$, $b=2\beta$, and $p=2q$, we obtain
\best
C=\left|{\cal L}^{123}\right|= \sum_{m=4\alpha \beta}\;
\sum_{0\leq 2q\leq 2\alpha-1} 1=\sum_{\alpha|{m\over 4}} \alpha
=\sigma\l({m\over 4}\r).\qquad \Box
\eest

\medskip

\begin{prop}
\label{rel RT-Ta}
The generating function $\phi_C$ of an embedded torus $C$ is
\bear
\phi_C(t)=\l[f(t)f(t^2)^{s_1/ 2}f(t^4)^{s_2/4}\r]^{\mbox{{\em sgn} C}}
\label{proprelRT-Ta}
\eear
where
\bear\label{s_i's}
s_1=\ma\sum_{i=1}^3 \sgn\,D_i-3, \qquad
s_2=\ma\prod_{i=1}^3 \sgn\,D_i-\ma\sum_{i=1}^3 \sgn\,D_i+2,
\eear
and
\bear
f(t)=\prod_{m\ge 1} F(t^m)^{\sigma(m)/m}.
\label{deff(t)}
\eear
\end{prop}
\pf From  equations (\ref{phiC}) and (\ref{RTm3}) and Lemma
\ref{abclemma}, we obtain
\best
\log \phi_C(t_A)=
\mbox{sgn}\; D_0 \ma \sum_{m= 1}^\infty
\frac1m\left[\sigma(m)+s_1\cdot\sigma\l({m\over 2}\r)+
s_2\cdot\sigma\l({m\over 4}\r)\right]\  \log F(t_A^m)
\eest
After substituting in (\ref{deff(t)}), this gives (\ref{proprelRT-Ta}).
\qed

\bigskip

When $C$ has Taubes' type 0 all three $D_i$ have positive sign, so
$s_1=s_2=0$ in (\ref{s_i's}).  Similarly,
$(s_1,s_2)$ is $(-2,0)$ for type 1, $(-4,4)$  for type 2,  and $(-6,4)$ for
type 3.  Thus
(\ref{proprelRT-Ta}) gives
$$
 \prod_{A\in T} \ \ma \prod_{C\in\M_A} \phi_C(t_A)\ =\  \prod_{A\in T} \
\ma\prod_{ k=0}^3 f_k(t_A)^{\tau(A,k)}
$$
where $f_0=f$ and
$$
f_1(t) ={f(t)\over f(t^2)},\qquad
f_2(t)={f(t)f(t^4)\over f^2(t^2)},\qquad
f_3(t)={f(t)f(t^4)\over f^3(t^2)}
$$
--- exactly as in (\ref{T gen fc})! Since the first factors in
(\ref{egfgeneratingfnc}) and (\ref{3.defgr}) are equal by  (\ref{4.1}),
we have the following equivalence.

\begin{theorem}
For any closed symplectic four-manifold $(X,\w)$, the Taubes and Ruan-Tian
Series (\ref{egfgeneratingfnc}) and
(\ref{3.defgr}) coincide when $f$ and $F$ are related by (\ref{deff(t)}):
$$
GT_X(t,s)\ =\ RT_X(t,s).
$$
Hence Taubes' Gromov invariants $GT_{\delta}(A)$ depend
only of the deformation class of
$\w$ and are computable from the Ruan-Tian invariants.
\label{4.mainthm}
\end{theorem}

\medskip

If we use the particular form of $F$ satisfying  (\ref{prodF=t})
and  make the change of variable
$m=ab$, we obtain
\best
\log f(t)&=&\sum_{m= 1}^\infty \;\sum_{a|m}{a\over m}\log F(t^m)
= \; \sum_{b= 1}^\infty {1\over b}\;\sum_{a=1}^\infty \log F(t^{ab})=
 \ma \sum_{b= 1}^\infty {1\over b}\; t^b=\;\log {1\over 1-t}.
\eest
Thus the Ruan-Tian Series with the generating functions defined in
(\ref{MoebiusF})
exactly reproduces the Taubes Series with his choice of generating
functions (\ref{2.taubesgeneratingfncs}).

\bigskip


\end{document}